\documentclass{amsart}
\usepackage{amssymb,amscd}

\newtheorem{lemma}{Lemma}[section]
\newtheorem{theorem}[lemma]{Theorem}
\newtheorem{proposition}[lemma]{Proposition}
\newtheorem{corollary}[lemma]{Corollary}

\newtheorem{problem}[lemma]{Problem}

\newcommand{\F}{\mathcal {F}}

\newcommand{\U}{\mathcal U}
\newcommand{\V}{\mathcal V}

\newcommand{\A}{\mathcal A}
\newcommand{\B}{\mathcal B}

\newcommand{\IN}{\mathbb N}
\newcommand{\IZ}{\mathbb Z}

\newcommand{\inv[1]}{\overset{\leftrightarrow}{#1}}

\newcommand{\LL}{\mathcal L}

\title[Algebra in the superextensions of groups, III]{Algebra in the superextensions of groups, III: minimal left ideals of $\lambda(\mathbb Z)$}
\author{T.Banakh, V.Gavrylkiv}
\address{Department of Mathematics, Ivan Franko National University of Lviv and Instytut Matematyki, Uniwersytet Humanistyczno-Przyrodniczy w Kielcach}
\email{tbanakh@yahoo.com}
\address{Vasyl Stefanyk Precarpathian National University, Ivano-Frankivsk (Ukraine)}
\email{vgavrylkiv@yahoo.com}

\begin{document}

\begin{abstract}
We prove that the minimal left ideals of the superextension
$\lambda(\IZ)$ of the discrete group $\IZ$ of integers are metrizable
topological semigroups, topologically isomorphic to minimal left ideals of the superextension $\lambda(\IZ_2)$ of the compact group $\IZ_2$ of integer 2-adic numbers. 

The superextension $\lambda(X)$ of a discrete group
$X$ is the compact Hausdorff right-topological semigroup
consisting of maximal linked systems on $X$ and endowed with the
semigroup operation $\A*\B=\{A\subset X:\{x\in
X:x^{-1}A\in\B\}\in\A\}.$
\end{abstract}
\maketitle
\section*{Introduction}

After the topological proof (see \cite[p.102]{HS}, \cite{H2}) of
Hindman theorem \cite{Hind}, topological methods become a standard
tool in the modern combinatorics of numbers, see \cite{HS},
\cite{P}. The crucial point is that any semigroup operation $\ast$
defined on any discrete space $X$ can be extended to a
right-topological semigroup operation on $\beta(X)$, the Stone-\v
Cech compactification of $X$. The extension of the operation from
$X$ to $\beta(X)$ can be defined by the simple formula:
\begin{equation}\label{extension}\U\ast\V=\Big\{\bigcup_{x\in
 U}x*V_x:U\in\U,\;\{V_x\}_{x\in U}\subset\V\Big\},
\end{equation}
where $\U,\V$ are ultrafilters on $X$. Endowed with the
so-extended operation, the Stone-\v Cech compactification
$\beta(X)$ becomes a compact right-topological semigroup. The
algebraic properties of this semigroup (for example, the existence
of idempotents or minimal left ideals) have important consequences
in combinatorics of numbers, see \cite{HS}, \cite{P}.

The Stone-\v Cech compactification $\beta(X)$ of $X$ is the
subspace of the double power-set $\mathcal P(\mathcal P(X))$,
which is a complete lattice with respect to the operations of
union and intersection. In \cite{G2} it was observed that the
semigroup operation extends not only to $\beta(X)$ but also to the
complete sublattice $G(X)$ of $\mathcal P(\mathcal P(X))$, 
generated by $\beta(X)$. This complete sublattice consists of all
inclusion hyperspaces over $X$.

By definition, a family $\F$ of non-empty subsets of a discrete
space $X$ is called an {\em inclusion hyperspace} if $\F$ is
monotone in the sense that a subset $A\subset X$ belongs to $\F$
provided $A$ contains some set $B\in\F$. Besides the operations of union
and intersection, the set $G(X)$ possesses an important transversality operation assigning to each
inclusion hyperspace $\F\in G(X)$ the inclusion hyperspace
$$\F^\perp=\{A\subset X:\forall F\in\F\;(A\cap F\ne\emptyset)\}.$$
This operation is involutive in the sense that $(\F^\perp)^\perp=\F$.

It is known that the family $G(X)$ of inclusion hyperspaces on $X$
is closed in the double power-set $\mathcal P(\mathcal
P(X))=\{0,1\}^{\mathcal P(X)}$ endowed with the natural product
topology. The induced topology on $G(X)$ can be described
directly: it is generated by the sub-base consisting of the sets
$$U^+=\{\F\in G(X):U\in\F\}\mbox{ and }U^-=\{\F\in G(X):U\in\F^\perp\}$$
where $U$ runs over subsets of $X$. Endowed with this topology,
$G(X)$ becomes a Hausdorff supercompact  space. The latter means
that each cover of $G(X)$ by the sub-basic sets has a 2-element
subcover.

The extension of a binary operation $\ast$ from $X$ to $G(X)$ can
be defined in the same way as for ultrafilters, i.e., by the
formula~(\ref{extension}) applied to any two inclusion hyperspaces
$\U,\V\in G(X)$. In \cite{G2} it was shown that for an associative
binary operation $\ast$ on $X$ the space $G(X)$ endowed with the
extended operation becomes a compact right-topological semigroup.
Besides the Stone-\v Cech extension, the semigroup $G(X)$ contains
many important spaces as closed subsemigroups. In particular,  the
space
$$\lambda(X)=\{\F\in G(X):\F=\F^\perp\}$$of maximal linked systems on
 $X$ is a closed subsemigroup of $G(X)$. The space $\lambda(X)$ is
  well-known in  General and Categorial Topology as the {\em
 superextension}
  of $X$, see \cite{vM}, \cite{TZ}.  Endowed with the extended binary
 operation, the
  superextension $\lambda(X)$ of a semigroup $X$ is a supercompact
  right-topological semigroup containing $\beta(X)$ as a subsemigroup.

The thorough study of algebraic properties of the superextensions
$\lambda(X)$ of  groups $X$ was started in \cite{BGN} and
continued in \cite{BG2}. In this paper we  concentrate at
describing the minimal (left) ideals of $\lambda(X)$.

Understanding the structure of minimal left ideals of the
semigroup $\beta(X)$ had important combinatorial consequences. For
example, properties of ultrafilters from a minimal left ideal of
$\beta(X)$ were exploited in the topological proof of the
classical Van der Waerden  Theorem \cite[14.3]{HS} due to
Fustenberg and Katznelson \cite{FK}. Minimal left ideals of the
semigroup $\beta(\IZ)$ play also an important role in Topological
Dynamics, see \cite{BB}, \cite{BF}, \cite[Ch.19]{HS}. We believe that studying the
structure of minimal (left) ideals of the semigroups $\lambda(X)$
also will have some combinatorial or dynamical consequences.


The main result of this paper is Theorem~\ref{minidealZ} asserting
that the minimal left ideals of the semigroup $\lambda(\IZ)$ are
compact metrizable topological semigroups topologically isomorphic to minimal
left ideals of the superextension $\lambda(\IZ_2)$ of the (compact
metrizable) group $\IZ_2$ of integer 2-adic numbers.

\section{Right-topological semigroups}

In this section we recall some information from \cite{HS} related
to right-topological semigroups. By definition, a
right-topological semigroup is a topological space $S$ endowed
with a semigroup operation $\ast:S\times S\to S$ such that for
every $a\in S$ the right shift $r_a:S\to S$, $r_a:x\mapsto x\ast
a$, is continuous. If the semigroup operation $\ast:S\times S\to
S$ is continuous, then $(S,\ast)$ is a {\em topological
semigroup}.


A non-empty subset $I$ of a semigroup $S$ is called a {\em left} (resp. {\em right})
{\em ideal\/} if $SI\subset I$ (resp. $IS\subset I$). If $I$ is both
a left and right ideal in $S$, then $I$ is called an {\em ideal}
in $S$. Observe that for every $x\in S$ the set $Sx=\{sx:s\in S\}$
(resp. $xS=\{xs:s\in S\}$) is a left (resp. right) ideal in $S$.
Such an ideal is called {\em principal}. An ideal $I\subset S$ is
called {\em minimal} if any ideal of $S$ that lies in $I$
coincides with $I$. By analogy we define minimal left and right
ideals of $S$. It is easy to see that each minimal left (resp.
right) ideal $I$ is principal. Moreover, $I=Sx$ (resp. $I=xS$) for
each $x\in I$.

If $S$ is a compact Hausdorff right-topological semigroup, then each
 minimal left ideal in $S$, being principal, is closed in $S$. By
\cite[2.6]{HS}, each left ideal in $S$ contains a minimal left ideal.
 The union of all minimal left ideals of $S$ coincides with the minimal
ideal $K(S)$ of $S$, \cite[2.8]{HS}. By \cite[2.11]{HS}, all the
minimal left ideals of $S$ are mutually homeomorphic.

An element $z$ of a semigroup $S$ is called a {\em right zero}
 in $S$ if $xz=z$  for all
$x\in S$. It is clear that $z\in S$ is a right zero in $S$
if and only if the singleton $\{z\}$ is a (minimal) left ideal in
$S$.

In the sequel we shall often use the following

\begin{lemma}\label{l1}

Let $X,Y$ be compact right-topological semigroups. If a semigroup
homomorphism $h:X\to Y$ is injective on some minimal left ideal of
$X$, then $h$ is injective on each minimal left ideal of $X$.
\end{lemma}

\begin{proof}
 Assume that $h$ is injective on a minimal left ideal $Xa$ of $X$ and take any
  other minimal left ideal $Xb$ of $X$. By \cite[2.11]{HS}, the right shift $r_a:X\to X$, $r_a:x\mapsto xa$,
  is injective on $Xb$. Next, consider the right shift $r_{h(a)}:Y\to Y$, $r_{h(a)}:y\mapsto y\cdot h(a)$.
  It follows from the equality $h\circ r_a=r_{h(a)}\circ h$ and the injectivity of the maps $r_a|Xb$
  and $h|Xa$  that the map $h|Xb$ is injective.
\end{proof}

\section{Inclusion hyperspaces and superextensions}

A family $\mathcal L$ of subsets of a set $X$ is called {\em a
linked system on $X$} if $A\cap B\ne\emptyset$ for all
$A,B\in\mathcal L$. Such a linked system $\mathcal L$ is {\em
maximal linked} if $\mathcal L$ coincides with any linked system
$\mathcal L'$ on $X$ that contains $\mathcal L$. Each
(ultra)filter on $X$ is a (maximal) linked system. A linked system
$\mathcal L$ on $X$ is maximal linked if and only if  for any
partition $X=A\cup B$ either $A$ or $B$ belongs to $\mathcal L$.

By $\lambda(X)$ we denote the family of all maximal linked systems
on $X$. Since each ultrafilter on $X$ is a maximal linked system,
$\lambda(X)$ contain the  Stone-\v Cech extension $\beta(X)$ of
$X$. It is easy to see that each maximal linked system on $X$ is
an inclusion hyperspace on $X$ and hence $\lambda(X)\subset G(X)$.
Moreover, it can be shown that $\lambda(X)=\{\A\in
G(X):\A=\A^\perp\}$. Let also $N_2(X)=\{\A\in
G(X):\A\subset\A^\perp\}$ denote the family of all linked
inclusion hyperspaces on $X$. By \cite{G1} both the subspaces
$\lambda(X)$ and $N_2(X)$ are closed in the compact Hausdorff
space $G(X)$.

Each function $f:X\to Y$ between sets $X,Y$ induces a continuous
map $Gf:G(X)\to G(Y)$ assigning to an inclusion
hyperspace $\A\in G(X)$ the inclusion hyperspace
$$Gf(\A)=\{B\subset Y:f^{-1}(B)\in\A\}\in G(Y).$$
The function $Gf$ maps $\lambda(X)$ into $\lambda(Y)$, so we can put $\lambda f=Gf|\lambda(X)$.

Given any semigroup operation $\ast:X\times X\to X$ on a set $X$
we can extend this operation to $G(X)$ letting
$$\U\ast\V=\Big\{\bigcup_{x\in U}x*V_x:U\in\U,\;\{V_x\}_{x\in U}\subset\V\Big\}$$
for inclusion hyperspaces $\U,\V\in G(X)$. Equivalently, the product $\U\ast\V$ can be defined as
$$\U\ast\V=\{A\subset X:\{x\in X:x^{-1}A\in \V\}\in\U\}$$
where $x^{-1}A=\{z\in X:x\ast z\in A\}$. By \cite{G2} the
so-extended operation turns $G(X)$ into a right-topological
semigroup. The structure of this semigroup was studied in details
in \cite{G2}. In particular, it was shown that for each group $X$
the minimal left ideals of $G(X)$ are singletons containing {\em
invariant} inclusion hyperspaces.

We call an inclusion hyperspace $\A\in G(X)$ {\em invariant} if
$x\A=\A$ for all $x\in X$. More generally, given a subgroup
$H\subset X$ we define $\A$ to be {\em $H$-invariant} if $x\A=\A$ for all
$x\in H$.

It follows from the definition of the topology on $G(X)$ that the
set $\inv[G](X)$ of invariant inclusin hyperspaces is closed in $G(X)$ and coincides with the  minimal ideal $K(G(X))$ of the semigroup $G(X)$.
Consequently, $K(G(X))$ is a closed rectangular
subsemigroup of $G(X)$. The {\em rectangularity} of $K(G(X))$
means that $\A\circ\mathcal B=\mathcal B$ for all $\A,\mathcal
B\in K(G(X))$.

\section{The minimal ideal of $\lambda(G)$ for odd groups}

In this section we characterize groups $G$ whose superextension $\lambda(G)$ has one-point minimal left ideals.

Following \cite{BGN}, we define a group $G$ to be {\em odd} if the order of each element $x$ of $G$ is odd.
If $G$ is a finite odd group, then the maximal linked system $$\LL=\{A\subset G:|A|>|G|/2\}$$ is invariant.
In fact, a group $G$ possesses an invariant maximal linked system if and only if $G$ is odd, see Theorem~3.2
of \cite{BGN}. By Proposition~3.1 of \cite{BGN}, a maximal linked system $\mathcal Z\in\lambda(G)$ on a group
 $G$ is invariant if and only if $\mathcal Z$ is a right zero of the semigroup $\lambda(G)$ if and only if the
 singleton $\{\mathcal Z\}$ is a minimal left ideal in $\lambda(G)$. Taking into account that the invariant
  maximal linked systems form a closed rectangular subsemigroup of $\lambda(G)$, we obtain the main result of this section.

\begin{theorem}\label{stideal}
A group $G$ is odd if and only if all the minimal left ideals of
$\lambda(G)$ are singletons. In this case the minimal ideal
$K(\lambda(G))$ of $\lambda(G)$ is a closed rectangular semigroup
consisting of invariant maximal linked systems.
\end{theorem}

Given a subgroup $H$ of a group $G$ let $G/H=\{xH:x\in G\}$ and
$\pi:G\to G/H$ denote the quotient map. It induces a continuous
map $\lambda\pi:\lambda(G)\to\lambda(G/H)$ between the
corresponding superextensions.

\begin{lemma}\label{inj}
For any $H$-invariant maximal linked system
$\A\in\lambda(H)\subset\lambda(G)$ the restriction of
$\lambda\pi:\lambda(G)\to\lambda(G/H)$ to the principal left ideal
$\lambda(G)*\A$ is injective.
\end{lemma}

\begin{proof}
Fix a section $s:G/H\to G$ of $\pi$. For every $\mathcal
L\in\lambda(G)$ let $\widetilde{\mathcal L}=\lambda\pi(\mathcal
L)\in\lambda(G/H)$ be the projection of $\LL$ onto $G/H$ and $\mathcal M=\lambda s(\widetilde{\mathcal L})\in\lambda(G)$ be the lift of $\widetilde\LL$ by the section $s$.

We claim that $\mathcal L*\A=\mathcal M*\A$. Since $\LL*\A$ and
$\mathcal M*\A$ are maximal linked systems, it suffices to check
that $\LL*\A\subset\mathcal M*\A$.  Take any set $\bigcup_{x\in
L}x*A_x\in\LL*\A$ where $L\in\LL$ and $\{A_x\}_{x\in L}\subset\A$.
Consider the set $M=s\circ\pi(L)\in\mathcal M$. For every point
$y\in M$ find a point $x_y\in L$ with $y=s\pi(x_y)$ and observe
that $yH=\pi(y)=\pi(x_y)=x_yH$, which implies $y^{-1}x_y\in H$ and
hence $y^{-1}x_yA_{x_y}\in\A$ by the $H$-invariantness of $\A$.
Since $$\mathcal M*\A\ni\bigcup_{y\in
M}y(y^{-1}x_y*A_{x_y})=\bigcup_{y\in M}x_y*A_{x_y}\subset
\bigcup_{x\in L}x*A_x$$we conclude that $\bigcup_{x\in
L}x*A_x\in\mathcal M*\A$.

Now we are able to prove that
$\lambda\pi:\lambda(G)*\A\to\lambda(G/H)$ is injective. Take any
two distinct elements $\LL_1*\A\ne\LL_2*\A$ of $\lambda(G)*\A$.
For every $i\in\{1,2\}$ consider the maximal linked systems
$\widetilde\LL_i=\lambda\pi(\LL_i)=\lambda\pi(\LL_i*\A)$ and
$\mathcal M_i=\lambda s(\widetilde\LL_i)$. It follows from
$\mathcal M_1*\A=\LL_1*\A\ne\LL_2*\A=\mathcal M_2*\A$ that
$\mathcal M_1\ne\mathcal M_2$ and hence
$$\lambda\pi(\LL_1*\A)=\widetilde
\LL_1\ne\widetilde\LL_2=\lambda\pi(\LL_2*\A).$$
\end{proof}

\begin{corollary}\label{injmap}
For a normal odd subgroup $H$ of a group $G$ the map
$\lambda\pi:\lambda(G)\to\lambda(G/H)$ is injective on each
minimal left ideal of $\lambda(G)$. Consequently, every minimal
left ideal of $\lambda(G)$ is topologically isomorphic to a
minimal left ideal of $\lambda(G/H)$.
\end{corollary}

\begin{proof}

By Lemma~\ref{l1}, it suffices to show that $\lambda\pi$ is
injective on some minimal left ideal. The group $H$, being odd,
admits an $H$-invariant maximal linked system
$\A\in\lambda(H)\subset\lambda(G)$. By Lemma~\ref{inj} the
homomorphism $\lambda\pi$ is injective on the left ideal
$\lambda(G)*\A$ and hence is injective on any minimal left ideal
contained in $\lambda(G)*\A$ (it exists because $\lambda(G)$ is a
compact right-topological semigroup).
\end{proof}

\section{Maximal invariant linked systems on groups}

As we have seen in the preceding section, the property of a
maximal system $\LL\in\lambda(G)$ to be invariant is very strong
and forces $\LL$ to be a right zero of $\lambda(G)$. Such maximal
linked systems exist only on odd groups.

On the other hand, maximal invariant linked systems exist on each
group. An invariant linked inclusion hyperspace
$\LL\in\inv[N]_2(G)$ is called a {\em maximal invariant linked
system} if $\LL=\LL'$ for any invariant linked inclusion
hyperspace $\LL'\in\inv[N]_2(G)$ enlarging $\LL$. By the Zorn
Lemma, each invariant linked inclusion hyperspace can be enlarged
to a maximal invariant linked system.

\begin{proposition}\label{p4.1}

For any maximal invariant linked system $\LL_0$ on a group $G$ the
set $${\uparrow}\LL_0=\{\LL\in\lambda(G):\LL\supset\LL_0\}$$ is a
left ideal in $\lambda(G)$.
\end{proposition}

\begin{proof}

Let $\A,\mathcal B\in\lambda(X)$ be maximal linked systems with
$\LL_0\subset \mathcal B$. Then for every subset $L\in\LL_0$ we
get $$L=\bigcup_{x\in G}x(x^{-1}L)\in\A*\LL_0\subset\A*\mathcal
B$$ which means that $\LL_0\subset\A*\mathcal B$.
\end{proof}

Observe that $\LL_0\subset\LL\subset\LL_0^\perp$ for every
$\LL\in{\uparrow}\LL_0$. The following theorem shows that the
difference $\LL_0^\perp\setminus\LL_0$ (and consequently,
$\LL\setminus\LL_0$) is relatively small (for the group $G=\IZ$ it
is countable!).

\begin{theorem}\label{narist}
If $\LL_0$ is a maximal invariant linked system on an Abelian
group $G$, then for any subset $A\in\LL_0^\perp\setminus\LL_0$
there is a point $x\in G$ such that $xA=G\setminus A$ and
consequently, $A=x^2A$.
\end{theorem}

\begin{proof} Fix a subset $A\in\LL_0^\perp\setminus\LL_0$.
We claim that
\begin{equation}\label{eq:a}
aA\cap A=\emptyset
\end{equation}
for some $a\in G$. Assuming the converse, we would conclude that
the family $\{xA:x\in G\}$ is linked and then the invariant linked
system $\LL_0\cup\{xA:x\in G\}$ is  strictly larger than $\LL_0$,
which impossible because of the maximality of $\LL_0$.

Next, we find $b\in G$ with
\begin{equation}\label{eq:b}A\cup bA=G.
\end{equation}
Assuming that no such a point $b$ exist, we conclude that for any
$x,y\in G$ the union $xA\cup yA\neq G$. Then $(G\setminus xA)\cap
(G\setminus yA)=G\setminus (xA\cup yA)\neq\emptyset$, which means
that the family $\{G\setminus xA: x\in G\}$ is linked and
invariant. We claim that $G\setminus A\in\LL_0^\perp$. Assuming
the converse, we would conclude that $G\setminus A$ misses some
set $L\in\LL_0$. Then $L\subset A$ and hence $A\in\LL_0$ which is
not the case. Thus $G\setminus A\in\LL_0^\perp$ and hence
$\{G\setminus xA:x\in G\}$ because $\LL_0^\perp$ is invariant.
Since $\LL_0\cup\{G\setminus xA:x\in G\}$ is an invariant linked
system containing $\LL_0$, the maximality of $\LL_0$ guarantees
that $G\setminus A\in\LL_0$ which contradicts $A\in\LL_0^\perp$.

Finally we show that $G\setminus A=aA=bA$. Observe that
(\ref{eq:a}) and (\ref{eq:b}) imply that $aA\subset bA$ and hence
$A\subset a^{-1}bA$. On the other hand, (\ref{eq:a}) and
(\ref{eq:b}) are equivalent to  $a^{-1}A\cap A=\emptyset$ and
$b^{-1}A\cup A=G$, which implies $a^{-1}A\subset b^{-1}A$ and this
yields $ba^{-1}A\subset A$. Unifying this inclusion with $A\subset
a^{-1}bA=ba^{-1}A$, we conclude that $ba^{-1}A=A$ and hence
$bA=aA$. Now looking at (\ref{eq:a}) and (\ref{eq:b}) we see that
$G\setminus A=aA=bA$.
\end{proof}

\section{Minimal left ideals of $\lambda(\IZ)$}
In this section we apply the results of the preceding sections to
describe the structure of minimal left ideals of the semigroup
$\lambda(\IZ)$. It turns out that they are isomorphic to minimal
left ideals of the superextension $\lambda(\IZ_2)$ of the compact
topological group $\IZ_2$ of integer 2-adic numbers. We recall
that $\IZ_2=\varprojlim C_{2^k}$ is a totally disconnected compact
metrizable Abelian group, which is the limit of the inverse
sequence
$$\dots\to C_{2^n}\to\dots\to C_8\to C_4\to C_2$$
of cyclic 2-groups $C_{2^n}$. Let $\pi:\IZ\to \IZ_2$ denote the
canonic (injective) homomorphism of $\IZ$ into $\IZ_2$ (induced by
the quotient maps $\pi_{2^k}:\IZ\to \IZ/2^k\IZ=C_{2^k}$, $k\in\IN$).

By the continuity of the functor $\lambda$ in the category of compact Hausdorff spaces (see \cite[2.3.2]{TZ}),
the superextension $\lambda(\IZ_2)$ can be identified with the
limit of the inverse sequence $$\dots\to
\lambda(C_{2^n})\to\dots\to \lambda(C_8)\to \lambda(C_4)\to
\lambda(C_2)$$ of finite semigroups $\lambda(C_{2^k})$. This
implies that $\lambda(\IZ_2)$ is a metrizable zero-dimensional
compact topological semigroup.

\begin{theorem}\label{minidealZ}
The homomorphism $\lambda\pi:\lambda(\IZ)\to\lambda(\IZ_2)$ is
injective on each minimal left ideal of $\lambda(\IZ)$.
Consequently, the minimal left ideals of the semigroup $\lambda(\IZ)$ are compact metrizable topological semigroups.
\end{theorem}

\begin{proof}
By Lemma~\ref{l1}, it suffices to check that the homomorphism
$\lambda\pi$ is injective on some minimal left ideal of
$\lambda(\IZ)$. Fix any maximal invariant linked system $\LL_0$ on
$\IZ$ (such a system exists by Zorn Lemma). By
Proposition~\ref{p4.1} the set
${\uparrow}\LL_0=\{\LL\in\lambda(\IZ):\LL\supset\LL_0\}$ is a left
ideal which necessarily contains a minimal left ideal $I$ of
$\lambda(\IZ)$. We claim that the homomorphism
$\lambda\pi:\lambda(\IZ)\to\lambda(\IZ_2)$ is injective on $I$.
Given two different maximal linked system $\A,\mathcal B\in I$ we
need to check that $\lambda\pi(\A)\ne\lambda\pi(\mathcal B)$.

Since the superextension $\lambda(\IZ_2)$ is the limit of the
inverse sequence $$\dots\to \lambda(C_{2^n})\to\dots\to
\lambda(C_8)\to \lambda(C_4)\to \lambda(C_2),$$ the inequality
$\lambda\pi(\A)\ne\lambda\pi(\mathcal B)$ will follow as soon as
we find $k\in\IN$ such that
$\lambda\pi_{2^k}(\A)\ne\lambda\pi_{2^k}(\B)$ where
$\lambda\pi_{2^k}:\lambda(\IZ)\to\lambda(C_{2^k})$ is the
homomorphism induced by the quotient homomorphism
$\pi_{2^k}:\IZ\to C_{2^k}$.

Pick any set $A\in\A\setminus\B$. Since $A\in
\LL_0^\perp\setminus\LL_0$, we can apply Theorem \ref{narist} to
conclude that $A=2n+A$ for some positive number $n\in \IZ$. The
later equality means that $A=\pi_{2n}^{-1}(\pi_{2n}(A))$ is the
complete preimage of the set $\pi_{2n}(A)$ under the quotient
homomorphism $\pi_{2n} : \IZ \to \IZ/2n\IZ=C_{2n}$. It follows
that
$\pi_{2n}(A)\in\lambda\pi_{2n}(\A)\setminus\lambda\pi_{2n}(\mathcal
B)$ and hence $\lambda\pi_{2n}(\A)\ne\lambda\pi_{2n}(\mathcal B)$.

Write the number $2n$ as the product $2n=2^k\cdot m$ for some odd
number $m$ and find a (unique) subgroup $H\subset C_{2n}$ of order
$|H|=m$. It follows that the quotient group $C_{2n}/H$ can be
identified with the cyclic 2-group $C_{2^k}$ so that  $q\circ
\pi_{2n}=\pi_{2^k}$ where $q:C_{2n}\to C_{2n}/H=C_{2^k}$ is the
quotient homomorphism. Corollary~\ref{injmap} guarantees that the
homomorphism $\lambda q:\lambda(C_{2n})\to\lambda(C_{2^k})$ is
injective on each minimal left ideal of $\lambda(C_{2n})$. In
particular, it is injective on the minimal left ideal
$\lambda\pi_{2n}(I)$. Consequently, $\lambda\pi_{2^k}(\A)=\lambda
q(\tilde \A)\ne\lambda q(\tilde\B)=\lambda\pi_{2^k}(\B)$. This
completes the proof of the injectivity of
$\lambda\pi:\lambda(\IZ)\to\lambda(\IZ_2)$ on the left ideal $I$
and consequently, on each minimal left ideal $J$ of
$\lambda(\IZ)$.

Since minimal left ideals of $\lambda(\IZ)$ are compact, the
restriction $\lambda\pi|J$ is a topological isomorphism of $J$
onto the minimal left ideal $\lambda\pi(J)$ of $\lambda(\IZ_2)$.
Since $\lambda(\IZ_2)$ is a metrizable topological semigroup, so
are the semigroups $\lambda\pi(J)$ and $J$.
\end{proof}

\section{Some Open Problems}

We saw in Theorem~\ref{stideal} that the minimal ideal
$K(\lambda(G))$ of the superextension of an odd group $G$ is a
compact topological semigroup.

\begin{problem}
Characterize groups $G$ such that the minimal ideal
$K(\lambda(G))$ is closed in $\lambda(G)$. Is the minimal ideal
$K(\lambda(\IZ))$ closed in $\lambda(\IZ)$? Is $K(\lambda(\IZ))$ a
topological semigroup?
\end{problem}

\begin{problem}
Characterize groups $G$ such that the minimal left ideals of
$\lambda(G)$ are (metrizable) topological semigroups.
\end{problem}

\end{document}